\documentclass[11pt,bezier]{article}
\usepackage{amsfonts}

\title{On the Upper Bounds of MDS Codes}
\author{Jiansheng Yang\footnote{Supported by Shanghai Leading Academic Discipline Project£¬
Project Number£ºS30104. Email: yjsyjs@staff.shu.edu.cn} ,\hspace{3mm} Yunying Zhang\\ \\
Department of Mathematics, Shanghai University \\
 Shanghai 200444, China }
\date{}
\begin{document}
\maketitle

\begin{abstract}
Let $M_{q}(k)$ be the maximum length of MDS codes with parameters
$q,k$. In this paper, the properties of $M_{q}(k)$ are studied, and
some new upper bounds of $M_{q}(k)$ are obtained. Especially we
obtain that $M_{q}(q-1)\leq q+2(q\equiv4(mod~6)),~M_{q}(q-2)\leq
q+1(q\equiv4(mod~6)),~M_{q}(k)\leq q+k-3~(q=36(5s+1),~s\in N$ and $
k=6,7).$

\end{abstract}

\section{Introduction}

\indent

Let $C$ be an $(n,q^{k},d)$ code, if $d=n-k+1$, then $C$ is called a
maximum distance separable (MDS) code. MDS codes are at the heart of
combinatorics and finite geometries. In their book [9] MacWilliams
and Sloane describe MDS codes as ¡°one of the most fascinating
chapters in all of coding theory.¡± These codes can be linear or
non-linear. Very little is known about non-linear $(n,q^{k},n-k+1)$
MDS codes. R.H.Bruck , H.J.Ryser, R.Silverman and A.A.Bruen had
proved some results on MDS codes [4,5,10] early. Recently,
T.L.Alderson studies MDS codes extension. And he has obtained some
important results[1,2,3].

In this paper, we assume that $A$=$\{0,1,2,\cdots,q-1\}$ is a
additive group(not necessary cyclic group). Denote the maximum
number $n$ of an $(n,q^{k},n-k+1)$ MDS code over $A$ by $M_{q}(k)$.
If $A$ is a field, then denote the maximum number $n$ of a linear
$(n,q^{k},n-k+1)$ MDS code over $A$ by $m_{q}(k)$. The Main
Conjecture of $m_{q}(k)$ is the following.
$$m_q (k) = \left\{
{\begin{array}{*{20}c}
   {q + 2} & {for\;k = 3\;and\;k = q - 1\;both\;q\;even,}  \\
   {q + 1} & {in\;all\;other\;cases.\;\quad \quad \quad \quad \quad \quad \quad }  \\
\end{array}} \right.$$

The Main Conjecture has been proved in some cases, [7] gave us a
good summarize. For $M_{q}(k)$, it is well known that $M_q(k)\leq
q+k-1$. In [5], A.A.Bruen and R.Silverman proved that:

\medskip

\textbf{Theorem 1.1} [5] (1) If $C$ is a $(q+k-1;k)$-MDS code with
$k\geq3$ and $q > 2$, then $4$ divides $q$.

(2) If $C$ is a $(q+k-1;k)$-MDS code with $k>3$ and $q>2$, then $36$
divides $q$.

\medskip

In [2], T.L.Alderson proved that:

\medskip

\textbf{Theorem 1.2} [2] (1) If 36 does not divide $q$ and $k\geq
4$, then a $q-ary~(n,k)$-MDS code satisfies $n\leq q+k-3$.

(2) If $q>2$ and $q\equiv 2~mod~4$ then no q-ary $(q+1; 3)$-MDS
codes exist.

\medskip

In [11], Wang proved that

\medskip

\textbf{Theorem 1.3} [11] $M_q(q-1)\leq q+1$ for $q$ is odd.

\medskip

In this paper, we use the generalized weight enumerator(see below)
and combinatorial methods to study $M_q(k)$. Some new upper bounds
of $M_{q}(k)$ are obtained. Especially we obtain that
$M_{q}(q-1)\leq q+2(q\equiv4(mod~6)),~M_{q}(q-2)\leq
q+1(q\equiv4(mod~6)),~M_{q}(k)\leq q+k-3~(q=36(5s+1),~s\in N$ and $
k=6,7).$

For the convenient, we introduce some notations and results as
following.

Let $C$ and $D$ be two codes of length $n$ over $A$. If there exist
$n$ permutations $\pi_{1},\cdots,\pi_{n}$ of the $q$ elements and a
permutation $\sigma$ of the $n$ coordinate positions such that
$(u_{1},\cdots,u_{n})\in
C$~iff~$(\pi_{1}(u_{\sigma{(1)}}),\cdots,\pi_{n}(u_{\sigma{(n)}}))\in
D$, then we call $C$ is equivalent to $D$. If $C$ is equivalent to
$D$, then $C$ and $D$ have the same Hamming distance. It is clear,
if $c_{0}\in A^{n}$, $D=c_{0}+C=\{c_{0}+\alpha|\alpha \in C\}$ is
equivalent to $C$. Then we may always assume that code $C$ contains
the zero element $\textbf{0}=(0,0,\cdots,0,0)$.

The generalized weight enumerator which is introduced by M.El-Khamy
and R.J.McEliece[6] is called the partition weight enumerator(PWE) .
Suppose the coordinate set $N=\{1,2,\cdots,n\}$ is partitioned into
$p$ disjoint subsets $N_{1},\cdots,N_{p}$, with $|N_{i}|=n_{i}$ ,
for $i=1,\cdots,p$. Denoting this partition by $\mathcal {T}$, the
$\mathcal {T}$-weight profile of an $v \in A^{n}$ is defined as
$\mathcal {W}_{\mathcal {T}}(v)=(\omega_{1},\cdots,\omega_{p})$,
where $\omega _{i}$ is the Hamming weight of $v$ restricted to
$N_{i}$. Given a code $C$ of length $n$, and an
$(n_{1},\cdots,n_{p})$ partition $\mathcal {T}$ of the $n$
coordinates of $C$, the $\mathcal {T}$-weight enumerator of $C$ is
defined as  following. $$A^{\mathcal
{T}}(\omega_{1},\cdots,\omega_{p})=|\{c \in C:\mathcal {W}_{\mathcal
{T}}(c)=(\omega_{1},\cdots,\omega_{p})\}|.$$

\medskip

{\bf Theorem 1.4} [6] For an $(n,q^{k},d)$ MDS code $C$ which
contains the zero element, the $p$-partition weight enumerator is
given by $$A^{\mathcal
{T}}(\omega_{1},\cdots,\omega_{p})=E(\omega)\frac{{n_{1} \choose
\omega_{1}}{n_{2} \choose \omega_{2}}\cdots{n_{p} \choose
\omega_{p}}}{{n \choose \omega}}$$\\where $\omega=\sum_{i=1}^n
\omega_{i}$, $E(\omega)=|\{c \in C:\mathcal {W}(c)=\omega\}|$ and
$\mathcal {W}(c)$ is the Hamming weight of $C$.

\medskip

{\bf Remark:} In the proof, [6] assume that $A$ is a field, however,
the proof is true for any $A$(It only need the MDS codes have the
zero element). Therefore, the formula holds for non-linear MDS codes
which contain the zero element.

\medskip

For an $(n,q^k,d)$ MDS code over $A$, the weight distribution is
known as $$E(\omega)=(q-1){n \choose
\omega}\sum_{j=0}^{\omega-d}(-1)^{j}{\omega-1 \choose
j}q^{\omega-d-j}$$ where $\omega\geq d$ [9], and we can know that
the formula holds for MDS codes(which contain the zero element) not
only for linear MDS codes[8].

For any $\alpha=(a_{1},\cdots,a_{n})\in A^{n}$, define the support
of $\alpha$ by $Supp\alpha=\{i|a_{i}\neq0,~1\leq i\leq n \}$, and
$\overline{Supp}\alpha=\{j|a_{j}=0,~1\leq j\leq n \}$.

\section{ New Upper Bounds for MDS Codes}

\indent

{\bf Theorem 2.1}~~If $q\equiv4(mod~6)$, then $M_{q}(q-1)\leq q+2$.

\medskip

\textbf{Proof:}~~Suppose $C$ is an $(q+3,q^{q-1},5)(q~is~even)$ MDS
code which contains the zero element. The
 partition $\mathcal {T}$ is given as following.$$\mathcal {T}=\mathcal {T}_{1}\cup\mathcal {T}_{2},~\mathcal {T}_{1}=\{1,2,3\},~\mathcal {T}_{2}=\{4,5,\cdots,q+3\}.$$
 $$n_{1}=|\mathcal {T}_{1}|=3,~n_{2}=|\mathcal {T}_{2}|=q,~\omega_{1}=2,~\omega_{2}=3.$$
 From Theorem 1.4, we have$$A^{\mathcal {T}}(2,3)=E(5)\frac{{3 \choose 2}{q \choose 3}}{{q+3 \choose 5}}=\frac{3q(q-1)^{2}(q-2)}{6}.$$
 where $E(5)=(q-1){q+3 \choose 5}$. For $(x,y)$, there are
 altogether $(q-1)^{2}$ pairs $(x,y)$ with $x,y\in S$ where
 $S=\{1,2,\cdots,q-1\}$. Thus there exists $(a,b)\in S$ such
 that$$|C_{a,b,0}|\geq\frac{3q(q-1)^{2}(q-2)}{6\times3(q-1)^{2}}=\frac{q(q-2)}{6}.$$
where

$$C_{a,b,0}=\{(a_{1},a_{2},a_{3},\cdots, a_{q+3})\in C|a_{1}=a, a_{2}=b,
a_{3}=0, a_{k}\in A, k=4,5,\cdots,q+3\}.$$

Since $C$ is an MDS code, w.l.g. we may assume $a=1,b=1$. Then
$$C_{1,1,0}=\{(1,1,0,a_{4},\cdots,a_{q+3})\in C|a_{k}\in A, k=4,5,\cdots,q+3\}.$$Let
$$C_{i}=C_{1,1,0,i}=\{(1,1,0,a_{4},\cdots,a_{q+3})\in C|a_{i}\neq0\}~(i\in\{4,5,\cdots,q+3\}).$$

We will prove that $|C_{1,1,0}|=\frac{q(q-2)}{6}$ when $q$ is even.
If this is not true, we have $ |C_{1,1,0}|> \frac{q(q-2)}{6}$, then
there exists $i$, w.l.g. assume $i$=4, such
that$$|C_{4}|\geq\frac{|C_{1,1,0}|{3 \choose 1}}{{q \choose
1}}>\frac{3q(q-2)}{6q}=\frac{q-2}{2}.$$

Assume
$\alpha=(1,1,0,a_{4},\cdots,a_{q+3}),~\beta=(1,1,0,b_{4},\cdots,b_{q+3})\in
C_{4}$, we have $a_{4}\neq0,~b_{4}\neq0$. If $i\in Supp\alpha \cap
Supp\beta~(5\leq i\leq q+3)$, then we have $d(\alpha,\beta)\leq4$, a
contradiction. Thus we have $Supp\alpha \cap Supp\beta=\{1,2,4\}$.
Let $\alpha_{1},\alpha_{2},\cdots,\alpha_{t}\in C_{4}$. Since
$Supp\alpha_{i}\cap Supp\alpha_{j}=\{1,2,4\}$ for all $i\neq j$ and
$\omega(\alpha)=5$, we have $\cup_{i=1}^t|Supp\alpha_{i}|=2t+3$.
This implies $2t+3\leq q+2$, i.e. $t\leq\frac{q-1}{2}$. Since $q$ is
even, we have $t\leq\frac{q-2}{2}$.

Hence$|C_{4}|\leq\frac{q-2}{2}~(q~is~even)$, a contradiction.

By this, we have$$|C_{1,1,0}|=\frac{q(q-2)}{6}~(q~is~even).$$ Thus
$\frac{q(q-2)}{6}$ must be an integer, however, if
$q\equiv4(mod~6)$, $\frac{q(q-2)}{6}$ is not an integer. Therefore,
if $q\equiv4(mod~6)$, then $M_{q}(q-1)\leq q+2$.\hfill$\Box$

\medskip

{\bf Theorem 2.2}~~If $q$ is even and $(l+2)!$ does not divide
$(q+l-1)\cdots(q+1)q(q-2)$ where $l\geq1$, then $M_{q}(q-2)\leq
q+l$.

\medskip

 \textbf{Proof:}~~Suppose $C$ is an $(q+l+1,q^{q-2},l+4)(q~is~even)$ MDS code which contains the zero element. The
 partition $\mathcal {T}$ is given as following.$$\mathcal {T}=\mathcal {T}_{1}\cup\mathcal {T}_{2},~\mathcal {T}_{1}=\{1,2\},~\mathcal {T}_{2}=\{3,4,\cdots,q+l+1\}.$$
 $$n_{1}=|\mathcal {T}_{1}|=2,~n_{2}=|\mathcal {T}_{2}|=q+l-1,~\omega_{1}=2,~\omega_{2}=l+2.$$
 From Theorem 1.4, we have$$A^{\mathcal {T}}(2,l+2)=E(l+4)\frac{{2 \choose 2}{q+l-1 \choose l+2}}{{q+l+1 \choose l+4}}={q+l-1 \choose l+2}(q-1).$$
 where $E(l+4)=(q-1){q+l+1 \choose l+4}$. For $(x,y)$, there are
 altogether $(q-1)^{2}$ pairs $(x,y)$ with $x,y\in S$ where
 $S=\{1,2,\cdots,q-1\}$. Thus there exists $(a,b)\in S$ such
 that$$|C_{a,b}|\geq\frac{{q+l-1 \choose l+2}(q-1)}{(q-1)^{2}}=\frac{{q+l-1 \choose l+2}}{q-1}.$$
where $C_{a,b}=\{(a_{1},a_{2},a_{3},\cdots,a_{q+l+1})|a_{1}=a,
a_{2}=b, a_{k}\in A, k=3,4,\cdots,q+l+1\}.$

 Since $C$ is an MDS code, we may assume $a=1,b=1$. Then
$$C_{1,1}=\{(1,1,a_{3},\cdots,a_{q+l+1})|a_{k}\in A, k=3,4,\cdots,q+l+1\}.$$Let
$B_{\underbrace{i,j,\cdots,k}_{l}}=C_{1,1,\underbrace{i,j,\cdots,k}_{l}}=\{(1,1,a_{3},\cdots,a_{q+l+1})|a_{k}\in
A, k=3,4,\cdots,q+l+1~and~a_{i},a_{j},\cdots,a_{k}\neq0\},$ where
$i,j,\cdots,k$ are the $l$ distinct numbers of
$\{3,4,\cdots,q+l+1\}.$

We claim that $|C_{1,1}|=\frac{{q+l-1 \choose l+2}}{q-1}$ when $q$
is even. If this is not true, since$ |C_{1,1}|\geq \frac{{q+l-1
\choose l+2}}{q-1}$, we have$|C_{1,1}|>\frac{{q+l-1 \choose
l+2}}{q-1}$ and there exist $i,j,\cdots,k$, w.l.g. assume the $l$
numbers are $3,4,\cdots,l+2$, such
that$$|B_{3,4,\cdots,l+2}|\geq\frac{|C_{1,1}|{l+2 \choose l}}{{q+l-1
\choose l}}>\frac{q-2}{2}.$$

Assume
$\alpha=(1,1,a_{3},\cdots,a_{q+l+1}),~\beta=(1,1,b_{3},\cdots,b_{q+l+1})\in
B_{3,4,\cdots,l+2}$, we have $a_{r}\neq0,~b_{s}\neq0,3\leq r,s\leq
l+2$. If $i\in Supp\alpha \cap Supp\beta~(l+3\leq i\leq q+l+1)$,
since $|Supp\alpha|=|Supp\beta|=l+4$, then we have
$d(\alpha,\beta)\leq l+3$, a contradiction. Thus we have $Supp\alpha
\cap Supp\beta=\{1,2,\cdots,l+2\}$. Let
$\alpha_{1},\alpha_{2},\cdots,\alpha_{t}\in B_{3,4,\cdots,l+2}$.
Since $Supp\alpha_{i}\cap Supp\alpha_{j}=\{1,2,\cdots,l+2\}$ for all
$i\neq j$ and $\omega(\alpha)=l+4$, we have
$\cup_{i=1}^t|Supp\alpha_{i}|=2t+l+2$. This implies $2t+l+2\leq
q+l+1$, i.e. $t\leq\frac{q-1}{2}$. Since $q$ is even, we have
$t\leq\frac{q-2}{2}$.

Hence$|B_{3,4,\cdots,l+2}|\leq\frac{q-2}{2}~(q~is~even)$, a
contradiction.

By this, we have$$|C_{1,1}|=\frac{{q+l-1 \choose
l+2}}{q-1}=\frac{(q+l-1)\cdots(q+1)q(q-2)}{(l+2)!}~(q~is~even).$$
Thus $\frac{(q+l-1)\cdots(q+1)q(q-2)}{(l+2)!}$ must be an integer.
Therefore, if $q$ is even and $(l+2)!$ does not divide
$(q+l-1)\cdots(q+1)q(q-2)$ where $l\geq1$, then $M_{q}(q-2)\leq
q+l$.\hfill$\Box$

\medskip

By calculating, we can get the following.

\medskip

Corollary 2.2.1 $M_{q}(q-2)\leq q+1~(q\equiv4(mod~6)).$

\medskip

Corollary 2.2.2 $M_{q}(q-2)\leq q+3~(q\equiv6~or~26(mod~30)).$

\medskip

Corollary 2.2.3 $M_{q}(q-2)\leq q+5~(q\equiv8~or~36(mod~42)).$

\medskip

{\bf Theorem 2.3}~~If $q$ is even and $(k-1)!$ does not divide
$(q+k-4)\cdots(q+1)q(q-2)$ where $k\geq4$, then $M_{q}(k)\leq
q+k-3$.

\medskip

 \textbf{Proof:}~~Suppose $C$ is an $(q+k-2,q^{k},q-1)(q~is~even)$ MDS code which contains the zero element. The
 partition $\mathcal {T}$ is given as following.$$\mathcal {T}=\mathcal {T}_{1}\cup\mathcal {T}_{2},~\mathcal {T}_{1}
 =\{1,2\},~\mathcal {T}_{2}=\{3,4,\cdots,q+k-2\}.$$
 $$n_{1}=|\mathcal {T}_{1}|=2,~n_{2}=|\mathcal {T}_{2}|=q+k-4,~\omega_{1}=2,~\omega_{2}=q-3.$$
 From Theorem 1.4, we have$$A^{\mathcal {T}}(2,q-3)=E(q-1)\frac{{2 \choose 2}{q+k-4 \choose q-3}}{{q+k-2 \choose q-1}}
 ={q+k-4 \choose q-3}(q-1).$$
 where $E(q-1)=(q-1){q+k-2 \choose q-1}$. For $(x,y)$, there are
 altogether $(q-1)^{2}$ pairs $(x,y)$ with $x,y\in S$ where
 $S=\{1,2,\cdots,q-1\}$. Thus there exists $(a,b)\in S$ such
 that\begin{equation}|C_{a,b}|\geq\frac{{q+k-4 \choose q-3}(q-1)}{(q-1)^{2}}=\frac{{q+k-4 \choose q-3}}{q-1}.\end{equation}
where $C_{a,b}=\{(a_{1},a_{2},a_{3},\cdots,a_{q+k-2})|a_{1}=a,
a_{2}=b, a_{m}\in A, m=3,4,\cdots,q+k-2\}.$

Since $C$ is an MDS code, we may assume $a=1,b=1$. Then
$$C_{1,1}=\{(1,1,a_{3},\cdots,a_{q+k-2})|a_{m}\in A, m=3,4,\cdots,q+k-2\}.$$Let
$B_{\underbrace{i,j,\cdots,r}_{k-3}}=C_{1,1,\underbrace{i,j,\cdots,r}_{k-3}}=\{(1,1,a_{3},\cdots,a_{q+k-2})|a_{m}\in
A,~m=3,4,\cdots,q+k-2~and~a_{i}=a_{j}=\cdots=a_{r}=0\},$ where
$i,j,\cdots,r$ are the $k-3$ distinct numbers of
$\{3,4,\cdots,q+k-2\}.$

We will prove that $|C_{1,1}|=\frac{{q+k-4 \choose q-3}}{q-1}$ when
$q$ is even. If this is not true, we have $|C_{1,1}|>\frac{{q+k-4
\choose q-3}}{q-1}$ and there exist $i,j,\cdots,r$, w.l.g. assume
the $k-3$ numbers are $3,4,\cdots,k-1$, such
that$$|B_{3,4,\cdots,k-1}|\geq\frac{|C_{1,1}|{k-1 \choose
k-3}}{{q+k-4 \choose k-3}}>\frac{q-2}{2}.$$

Assume
$\alpha=(1,1,a_{3},\cdots,a_{q+k-2}),~\beta=(1,1,b_{3},\cdots,b_{q+k-2})\in
B_{3,4,\cdots,k-1}$, we have $a_{r}=0,~b_{s}=0,3\leq r,s\leq k-1$.
If $i\in \overline{Supp}\alpha \cap \overline{Supp}\beta~(k\leq
i\leq q+k-2)$, then we have $d(\alpha,\beta)\leq q-2$, a
contradiction. Thus we have $\overline{Supp}\alpha \cap
\overline{Supp}\beta=\{3,4,\cdots,k-1\}$. Let
$\alpha_{1},\alpha_{2},\cdots,\alpha_{t}\in B_{3,4,\cdots,k-1}$.
Since $\overline{Supp}\alpha_{i} \cap
\overline{Supp}\alpha_{j}=\{3,4,\cdots,k-1\}$ for all $i\neq j$ and
$\omega(\alpha)=q-1$, we have $\cup_{i=1}^t|Supp\alpha_{i}|=2t+k-3$.
This implies $2t+k-3\leq q+k-4$, i.e. $t\leq\frac{q-1}{2}$. Since
$q$ is even, we have $t\leq\frac{q-2}{2}$.

Hence$|B_{3,4,\cdots,k-1}|\leq\frac{q-2}{2}~(q~is~even)$, a
contradiction.

By this,we have$$|C_{1,1}|=\frac{{q+k-4 \choose
q-3}}{q-1}=\frac{(q+k-4)\cdots(q+1)q(q-2)}{(k-1)!}~(q~is~even)$$
Thus $\frac{(q+k-4)\cdots(q+1)q(q-2)}{(k-1)!}$ must be an integer.
Therefore, if $q$ is even and $(k-1)!$ does not divide
$(q+k-4)\cdots(q+1)q(q-2)$ where $k\geq4$, then $M_{q}(k)\leq
q+k-3$.\hfill$\Box$

\medskip

By the Theorem 2.3, we have the following.

\medskip

Corollary 2.3.1 $M_{q}(k)\leq q+k-3~(q=36(5s+1),~s\in N$ and $
k=6,7).$

\section{Conclusion}

In this paper, we use the generalized weight enumerator and
combinatorial methods to study $M_q(k)$ which denote the maximum
number $n$ of an $(n,q^{k},n-k+1)$ MDS code. Compared to Theorem
1.1, Theorem 1.2 and Theorem 1.3, we obtain some new upper bounds of
$M_{q}(k)$. Especially we obtain that $M_{q}(q-1)\leq
q+2(q\equiv4(mod~6)),~M_{q}(q-2)\leq
q+1(q\equiv4(mod~6)),~M_{q}(k)\leq q+k-3~(q=36(5s+1),~s\in N$ and $
k=6,7).$

\end{document}